\documentclass[11pt,fleqn]{article}
\RequirePackage[T1]{fontenc}
\RequirePackage{amsthm,amsmath}
\RequirePackage[square,numbers]{natbib}
\usepackage[english]{babel}
\RequirePackage[colorlinks,citecolor=blue,urlcolor=blue]{hyperref}
\usepackage{enumerate}
\usepackage{amsmath,amsfonts,amssymb,euscript,mathrsfs}
\usepackage{fullpage}
\usepackage{amsmath}
\usepackage{dsfont}
\makeatletter
\def\captionof#1#2{{\def\@captype{#1}#2}}
\makeatother
\def\1{\mbox{\bf 1}}
\def\R{\mathbb{R}}

\def\N{\mathbb{N}}
\def\P{\mathbb{P}}
\def\E{\mathbb{E}}

\def\R{\mathbb{R}}

\def\Z{\mathbb{Z}}

\newtheorem{theo}{Theorem}
\newtheorem{lem}{Lemma}
\newtheorem{prop}{Proposition}

\newtheorem{Def/Prop}{Definition-Proposition}

\newcounter{exos}
\renewcommand\theexos{\arabic{exos}}

\newcounter{prob}
\renewcommand\theprob{\arabic{prob}}

\begin{document}
\title{Stationarity and ergodic properties for some observation-driven models in random environments}
\date{}

\author{Paul Doukhan \footnote{Universit\'e Cergy-Pontoise, UMR 8088 Analyse, G\'eom\'etrie et Mod\'elisation, 2 avenue Adolphe Chauvin, 95302 Cergy-Pontoise Cedex France.} 
\and
Michael H. Neumann \footnote{Friedrich-Schiller-Universit\"at Jena, Institut f\"ur Mathematik, Ernst-Abbe-Platz 2, 07743 Jena, Germany.}
\and
Lionel Truquet \footnote{UMR 9194 CNRS CREST, ENSAI, Campus de Ker-Lann, rue Blaise Pascal, BP 37203, 35172 Bruz cedex, France.}
 }

\maketitle

\begin{abstract}
The first motivation of this paper is to study stationarity and ergodic properties for a general class of time series models 
defined conditional on an exogenous covariates process. The dynamic of these models is given by an autoregressive latent process which forms a Markov chain in random environments. Contrarily to existing contributions in the field of Markov chains in random environments, the state space is not discrete and we do not use small set type assumptions or uniform contraction conditions for the random Markov kernels. 
Our assumptions are quite general and allow us to deal with models that are not fully contractive, such as threshold autoregressive processes. Using a coupling approach, we study the existence of a limit, in Wasserstein metric, for the backward iterations of the chain. We also derive ergodic properties for the corresponding skew-product Markov chain. Our results are illustrated with many examples
of autoregressive processes widely used in statistics or in econometrics, including  GARCH type processes, count autoregressions and categorical time series.     
\end{abstract}

\section{Introduction}

Non-linear time series have many important applications in various fields such as finance (\citep{tsay2}, \citep{mills}), economics (\citep{granger}) or climate analysis and ecology (\citep{silva}) among others.
Many textbooks now provide a thorough study of theoretical properties of non-linear autoregressive processes. See for instance \citep{Meyn}, \citep{D18} or \citep{Douc}.
A particularly tricky problem often encountered in studying non-linear time series models is to prove the existence of a stationary and ergodic path which is often the minimal condition needed for considering statistical applications such as likelihood inference.
Deriving stability properties of non-linear autoregressive processes using either Markov chain techniques or convergence results for iterated random systems has then attracted an important effort in the time series literature.

 However, it is difficult to find mathematical results concerning inclusion of covariates and in particular existence of stationary and ergodic paths when exogenous covariates are incorporated in the dynamic. In contrast, the use of exogenous covariates is almost systematic for practitioners who use such models.
In the applied statistical literature or in econometrics, there is a recent and growing interest in studying 
standard time series models with exogenous regressors and a few recent contributions already discussed this problem. 
See for instance \citep{GG}, \citep{Cav}, \citep{Francq}, \citep{deJong}, \citep{FT}, \citep{Truquet} or \citep{Debaly}. 
However, these contributions consider quite specific models and mainly with strong contraction assumptions in the sense that the transition kernels satisfy Lipschitz type properties with Lipschitz coefficients not depending on the exogenous process. 

The aim of this paper is to fill an important gap for this problem by deriving ergodic properties for more general non-linear structures
and only using weak contraction conditions. We focus on an important class of models called observation-driven.
An observation driven-model (of order $1$) is a bivariate time-homogeneous Markov chain $\left(\left(Y_t,\lambda_t\right)\right)_{t\in \Z}$ taking values in $\R^2$ or possibly in more general Cartesian products and where $\left(\lambda_t\right)_{t\in\Z}$ is an unobserved latent process such that for a Borel set $A$,
$$\P\left(Y_t\in A\vert \lambda_t\right)=p\left(A\vert \lambda_t\right),\quad \lambda_t=f\left(\lambda_{t-1},Y_{t-1}\right),$$
with $p\left(\cdot\vert\cdot\right)$ a probability kernel and $f$ a measurable function.
This important class of models is widely popular as it contains the well-known GARCH processes (\citep{Boug}, \citep{FZ}), models for time series of counts (see for instance \citep{FTR}, \citep{N}, \citep{DDM}, \citep{FTj}, \citep{DN}, \citep{Davis})
or models for categorical time series (see \citep{Fok}, \citep{FT}, \citep{Truquet}).
Deriving conditions under which there exists a unique invariant probability measure for such Markov chain models 
has then attracted a particular attention. An important difficulty arises for such discrete-valued process $(Y_t)_{t\in\Z}$ because the latent process $\left(\lambda_t\right)_{t\in\Z}$, which is itself a Markov chains with a non-discrete state space, does not satisfy the standard $\phi-$irreducibility properties. This leads various authors (\citep{FTR}, \citep{N}, \citep{DDM}, \citep{DN}) to develop new elegant methods to study these models, in particular coupling techniques or perturbation methods than can be applied to count autoregressions.

In this paper, we consider observation-driven models defined conditionally on a covariate process $X:=\left(X_t\right)_{t\in\Z}$. 
See Section \ref{S1} for a precise definition. In Econometrics, defining a dynamic conditional on an external stochastic process refers to the notion of strict exogeneity and the process $X$ is said to be strictly exogenous; see \citep{Sims} and \citep{Chamb} for a discussion about this notion. 
In probability, the notion of Markov chain in random environments is more likely used.
Unfortunately, the literature of Markov chains in random environments is not relevant for solving our problem. The seminal paper of \citep{Cogburn} only discusses discrete state spaces. The case of continuous state spaces is considered in \citep{Kifer} but the results are mainly appropriate for bounded state spaces and models that satisfy Doeblin's type condition, which is not adapted here.
Unbounded state spaces are considered in \citep{Stenflo} using contraction methods but the contraction condition is uniform with respect to $X$, a condition we want to relax. Moreover, we will also consider threshold models that are not fully contractive but semi-contractive. By semi-contractive, we mean that the function $f$ given above satisfies contraction properties only with respect to its first argument. This situation was considered recently in \citep{DN} and we will provide a non-trivial extension of their proof technique for studying observation-driven models in random environments. Let us mention that our approach could also be used to study higher-order observation-driven models
(i.e.~with several lag variables in the function $f$) but for readability, we prefer to only focus on first-order models which are the mostly used in practice.
 
The paper is organized as follow. In Section \ref{S1}, we give our main result. Its proof is provided in Section \ref{S2}.
In particular, the main difficulty is to study the convergence of the so-called backward iterations of the chain $\left(\lambda_t\right)_{t\in\Z}$ conditional on $X$ in Wasserstein metric. A discussion of the proof strategy, which is based on coupling, is given at the beginning of the section. Finally, we provide in Section \ref{S3} many examples of discrete and non-discrete observation-driven models satisfying our assumptions.  

\section{General result}\label{S1}
Let $g,k,d$ be three positive integers and $E,F,L$ be some Borel subsets of respectively $\R^g$, $\R^k$ and $\R^d$. We denote by $\mathcal{B}(E)$, $\mathcal{B}(F)$ and $\mathcal{B}(L)$ 
their corresponding Borel sigma-fields.
We consider a probability kernel $p$ from  $\left(F,\mathcal{B}(F)\right)$ to a $\left(E,\mathcal{B}(E)\right)$ as well as a stationary and ergodic stochastic process $X=(X_t)_{t\in\Z}$ taking values in $L$.
Our aim is to construct a process $\left((Y_t,\lambda_t)\right)_{t\in\Z}$ taking values in $E\times F$ and such that for $A\in \mathcal{B}\left(E\right)$ and $t\in\Z$,
\begin{equation}\label{dynamic}
\P\left(Y_t\in A\vert X,Y_{t-1},\lambda_t,Y_{t-2},\lambda_{t-1}\ldots\right)=p(A\vert \lambda_t),\quad\lambda_t=f\left(\lambda_{t-1},Y_{t-1},X_{t-1}\right),
\end{equation}
where $f$ is a measurable mapping from $F\times E\times L$ to $F$.

With such a formulation, the process $\left(\lambda_t\right)_{t\in\Z}$ is, conditionally on $X$, a time-inhomogeneous Markov chain.
In particular, if $h\colon\,F\rightarrow \R$ is a continuous and bounded function, we have
$$\E\left[h\left(\lambda_t\right)\vert \lambda_{t-1},X\right]=\int h\left(f(\lambda_{t-1},y,X_{t-1}\right)p\left(dy\vert \lambda_{t-1}\right).$$
We then introduce the random kernels $P_{X_t(\omega)}$ for $\omega\in\Omega$ and $t\in Z$ such that  
$$P_{X_t(\omega)}h(s)=\int h\left(f\left(s,y,X_t(\omega)\right)\right)p(dy\vert s).$$
Since our main goal is to study ergodic properties for the process $(Y_t)_{t\in\Z}$, we will also consider the transition kernel $R_{X_t}$, $t\in\Z$,
for the bivariate Markov chain in random environments $\left(Y_t,\lambda_t\right)_{t\in\Z}\vert X$. These kernels are defined by
$$R_{X_t(\omega)}h(y,s)=\int h\left(y', f\left(s,y,X_t(\omega)\right)\right)p\left(dy'\vert f\left(s,y,X_t(\omega)\right)\right)$$
for any continuous and bounded function $h\colon\,E\times F\rightarrow \R$. 

We remind that for two probability measures $\mu$ and $\nu$ defined on the same measurable space $\left(G,\mathcal{G}\right)$, their total variation distance is defined by 
$$d_{TV}\left(\mu,\nu\right)=\sup_{A\in \mathcal{G}}\left\vert\mu(A)-\nu(A)\right\vert.$$
We remind that the total variation distance can also be expressed in term of coupling,
$$d_{TV}\left(\mu,\nu\right)=\inf \P\left(X\neq Y\right),$$
where the infimum is taken over all pairs of random variables $(X,Y)$ such that $\P_X=\mu$ and $\P_Y=\nu$.
We also remind that if $P$ is a Markov kernel on $G$ and $\ell\colon\,G\rightarrow \R_+$ is a measurable function, the function $P\ell$ is defined by the equality $P\ell(z)=\int P(z,dz')\ell(z')$, $z\in G$.
Finally, for a nonnegative real number $u$, we set $\log^{+}(u)=\log\left(u\vee 1\right)$ where $u\vee 1=\max(u,1)$.
We will use the following assumptions. 

\begin{description}
\item[A1]
There exists a norm $\vert\cdot\vert$ on $F$ and a measurable function $\kappa\colon\,L\rightarrow \R_+$ such that $\E \log^{+}\kappa(X_0)<\infty$, $\E\log\kappa(X_0)<0$ and
that for all $y\in E$, $s,s'\in F$, $x\in L$ and $t\in\Z$,
$$\left\vert f(s,y,x)-f(s',y,x)\right\vert\leq \kappa(x)\vert s-s'\vert.$$

\item[A2] There exist three measurable functions $\gamma,\delta\colon\, L\rightarrow \R_+$, $V\colon\,L\rightarrow \R_+$ and $\alpha\in (0,1]$ such that $\E \log^{+}\delta(X_0)<\infty$, $\E\log^{+}\gamma(X_0)<\infty$, $\E\log \gamma(X_0)<0$, $V(s)\geq \vert s\vert^{\alpha}$ for $s\in L$ and
$$P_x V(s)\leq \gamma(x) V(s)+\delta(x)\mbox{ for } (s,x)\in F\times L.$$

\item[A3] There exists a polynomial function $\phi$, with positive coefficients, vanishing at $0$ and such that for every $(s,s')\in F^2$,
$$d_{TV}\left(p(\cdot\vert s),p(\cdot\vert s')\right)\leq 1- \exp\left(-\phi\left(\vert s-s'\vert\right)\right).$$
\end{description}

We now present our main result.

\begin{theo}\label{main}
Suppose that Assumptions {\bf A1-A3} hold true. Then there exists a stationary and ergodic process  $\left((Y_t,\lambda_t,X_t)\right)_{t\in\Z}$ solution of (\ref{dynamic}) and the probability distribution of such a process is unique.
\end{theo}

\paragraph{Discussion of the assumptions}
\begin{enumerate}
\item
Assumption {\bf A1} requires a contraction property for the function $f$ only with respect to its first argument. 
In particular, the function $f$ is not required to be continuous with respect to its second argument.
This semi-contractivity property will be particularly useful for defining threshold models or for getting sharp result for 
autoregressive categorical time series. For deterministic environments (i.e. without exogenous regressors), \citep{DN} recently used such a condition. We provide here an analogue for random environments.
Note that the condition $\E\log\left(\kappa(X_0)\right)<0$ is quite weak  
and cannot be removed in general for studying stationarity for such models. For instance, for the simple case $f(s,y,x)=a(x)s+b(x)$ where the $X_t'$s are i.i.d., the condition $\E\log\left(\kappa(X_0)\right)=\E\log\left\vert a(X_0)\right\vert<0$ is necessary for getting the existence of a stationary solution $\left(\lambda_t\right)_{t\in\Z}$; see \citep{BP}, Theorem $5$.

\item
Assumption {\bf A2} requires a drift condition for each Markov kernel $P_{X_t(\omega)}$. The drift parameters can be random 
and are allowed to only have a logarithmic moment. Condition $\E\log\left(\gamma(X_0)\right)<0$ is necessary and sufficient 
to ensure that products of type $\prod_{i=0}^n \gamma\left(X_{t-i}\right)$ vanish at infinity a.s. The latter property is important 
to get a standard drift condition of the form 
$$P_{X_{t-n}(\omega)}\cdots P_{X_t(\omega)}V\leq \eta\left(X_{t-n}(\omega),\ldots,X_t(\omega)\right)V+\zeta\left(X_{t-n}(\omega),\ldots,X_t(\omega)\right),$$
with $\eta\left(X_{t-n}(\omega),\ldots,X_t(\omega)\right):=\prod_{i=0}^n \gamma\left(X_{t-i}(\omega)\right)<1$ after iterating $n=n(\omega)$ successive random Markov kernels. 
Such a condition will be also central in our proof.

\item
Our proof of Theorem \ref{main} is based on the maximal coupling. In particular, a crucial step in our proof is to evaluate how close are two distributions $p\left(\cdot\vert s\right)$ and $p\left(\cdot\vert s'\right)$ in total variation. Assumption {\bf A3} requires these two distributions to be non singular unless the distance $\vert s-s'\vert$ goes to infinity.
For observation-driven models with deterministic environments, this kind of assumption has already been used in \citep{DN}, with a polynomial function $\phi$ of order one (i.e. $\phi(h)=K h$, $h\in\R_+$) is used. However, for some of the examples considered in the present paper, such as probit autoregressive processes or autoregressive processes with a Gaussian noise, the two distributions have to be asymptotically singular at a faster rate than the exponential rate. This explains the slightly more general assumption we use here.   
\end{enumerate}

\section{Proof of the main result}\label{S2}

\subsection{Maximal coupling and proof strategy}
Our approach will consist in studying first the Markov chain in random environments associated to the random Markov kernels
$P_{X_t(\omega)}$, $t\in \Z$.
In the spirit of the approach already used by \citet{Kifer} or \citet{Stenflo}, our aim is to show that under Assumptions {\bf A1-A3}, one can use a path-by-path approach and show
that for $\P-$almost all $\omega\in \Omega$, there exists a limit for the backward iterations $\delta_sP_{X_{t-n}(\omega)}\cdots P_{X_{t-1}(\omega)}$ when $n\rightarrow \infty$. To this end, we will use a Wasserstein metric. This limit, denoted by $\pi_t(\omega)$, will be shown to be a probability distribution on $\left(L,\mathcal{B}(L)\right)$ not depending on $s\in L$. Related to the dynamic (\ref{dynamic}), $\pi_t$ plays the role
of the conditional distribution of $\lambda_t$ given the covariate process $X$. It satisfies the invariance equation 
$\pi_tP_{X_t}=P_{X_{t+1}}$ a.s. 

To get such a result, we will use the maximal coupling and work first in the forward sense.  We define two processes $\left(\left(Y_t,\lambda_t\right)\right)_{t\geq 0}$ and $\left(\left(Y_t',\lambda'_t\right)\right)_{t\geq 0}$ and a probability measure $\overline{\P}_{\omega}$ such that $\lambda_0=s$, $\lambda'_0=s'$ and for $t\geq 0$,
$$\overline{\P}_{\omega}\left(Y_t\neq Y_t'\vert \lambda_t,\lambda_t'\right)=d_{TV}\left[p\left(\cdot\vert \lambda_t\right),p\left(\cdot\vert \lambda'_t\right)\right].$$ See for instance \citep{den}, Theorem $2.12$, 
for a proof of the existence of such a coupling.
We then define 
$$\lambda_{t+1}=f\left(\lambda_t,Y_{t-1},X_{t-1}(\omega)\right),\quad \lambda'_{t+1}=f\left(\lambda'_t,Y'_{t-1},X_{t-1}(\omega)\right).$$
We will also denote by $\overline{\E}_{\omega}$ the mathematical expectation corresponding to $\overline{\P}_{\omega}$.
Our first aim is to show that the probability
\begin{equation}\label{aborne} 
\overline{\P}_{\omega}\left(Y_{n+i}=Y'_{n+i}; i\geq 0 \mbox{ and } \sum_{i=1}^{\infty}\left\vert \lambda_{n+i}-\lambda'_{n+i}\right\vert \leq c_n(\omega)\right)
\end{equation}
is close to $1$ for a suitably chosen sequence $c_n(\omega)$ decreasing to $0$.

To this end, we will follow the approach used in \citet{DN} for deterministic environments. The important difficulty will be to adapt this approach taking into account the random parameters in Assumptions {\bf A1-A2}. 
The rest of the section will be organized as follows.

\begin{enumerate}
\item
In Section \ref{first}, we consider a subsampling of the bivariate time-inhomogeneous Markov chain $\left(\left(\lambda_t,\lambda'_t\right)\right)_{t\geq 0}$. In particular we introduce a sequence of some random times $\left(\tau_i(\omega)\right)_{i\geq 1}$ only depending on the random environment. Along this sequence, the drift parameters of this new Markov chain remain under control, as they are deterministic. We then control the tail of the distribution of some $d-$delayed return times $\rho_{\omega,j}$ (with $\rho_{\omega,j+1}-\rho_{\omega,j}\geq d=d_n(\omega)$) of the process $Z_{\omega}:=\left(\left(\lambda_{\tau_i(\omega)},\lambda'_{\tau_i(\omega)}\right)\right)_{i\geq 0}$ near the origin (more precisely in a ball with a non random radius). 
This will be obtained in Lemma \ref{annex}.
\item
When the process $Z_{\omega}$ is inside this ball at a given time $\rho_{\omega,j}$, we use Assumptions {\bf A1} and {\bf A3} 
to get a lower bound for the $\overline{\P}_{\omega}$-probability of the event $A_j:=\left\{Y_t=Y'_t\colon\, t\geq \tau_{\rho_{\omega,j}}(\omega)\right\}$. The semi-contraction condition {\bf A1} will be here of major importance. 
Then the probability that none of these events $A_j$ occur before time $n$ will be small when $n\rightarrow \infty$.
Moreover if such an event $A_j$ occurs, the difference $\left\vert \lambda_{\tau_{\rho_{\omega,j}}(\omega)}-\lambda_{\tau_{\rho_{\omega,j}}(\omega)}\right\vert$ will be quite small and this will help to fix the delays $d_n(\omega)$ and the rate $c_n(\omega)$ to control (\ref{aborne}). This will done in Section \ref{second}. 
\item
Finally, in Section \ref{third}, we use the lower bound obtained for the probability (\ref{aborne}) to control the convergence of the backward iterations of the chain and we prove Theorem \ref{main}.
\end{enumerate}

\subsection{Subsampling of the chain}\label{first}
For simplicity, we write $\kappa_t,\gamma_t,\delta_t$ for respectively $\kappa(X_t),\gamma(X_t),\delta(X_t)$.
\begin{lem}\label{control}
Assume that Assumptions {\bf A1-A3} hold true. Set for a positive integer $h$ and $t\in \Z$,
$$W_t^{(1)}=\delta_{t-1}+\sum_{i=1}^{\infty}\gamma_{t-1}\cdots\gamma_{t-i}\delta_{t-i-1},$$
$$W_t^{(2)}=\sup_{j\geq h}\gamma_{t-1}\cdots \gamma_{t-j},$$
$$W_t^{(3)}=\sup_{j\geq h}\kappa_{t-1}\cdots \kappa_{t-j},$$
$$W_t^{(4)}=\sum_{s=0}^{\infty}\phi\left(\kappa_{t+s}\cdots\kappa_t\right).$$
Then the process $\left(\left(W_t^{(1)},W_t^{(2)},W_t^{(3)},W_t^{(4)}\right)\right)_{t\in\Z}$ is stationary and ergodic. Moreover, if $h$ is large enough, there exists $C>1$ such that
$$\P\left(W_1^{(1)}\leq C, W_1^{(2)}\leq 1-1/C,W_1^{(3)}\leq 1-1/C,W_1^{(4)}\leq C\right)>0.$$
\end{lem}

\paragraph{Proof of Lemma \ref{control}}
Note first that stationarity and ergodic properties for the random vectors process follow from the representation 
$$G_t:=\left(W_t^{(1)},W_t^{(2)},W_t^{(3)},W_t^{(4)}\right)=H\left(\left(X_{t-j}\right)_{j\in\Z}\right),$$
for a suitable measurable function $H\colon\; F^{\Z}\rightarrow \R^4$. Stationarity and ergodicity of the process $\left(G_t\right)_{t\in\Z}$ then follows from that of $(X_t)_{t\in\Z}$.

Next, from the log-moment assumptions given in {\bf A2-A3}, 
the ergodic theorem ensures that 
$$\lim_{h\rightarrow \infty}\sup_{j\geq h}\gamma_{t-1}\cdots \gamma_{t-j}=0\quad\mbox{ a.s.}$$
A precise justification of this almost sure convergence can be found in \citep{Brandt}, see the proof of Lemma $1.1$.
The same property holds true if we replace $\gamma$ with $\kappa$.
Hence for $h$ sufficiently large, we have
$\P\left(W_1^{(2)}<1,W_1^{(3)}<1\right)>0$.
Next, we show that $W_t^{(1)}$ is finite a.s. Under the log-moments assumptions given in {\bf A2} and the stationarity and ergodicity of the process $\left(\left(\gamma_t,\delta_t\right)\right)_{t\in\Z}$, it is widely known that the stochastic recursions 
$$U_t=\gamma_{t-1}U_{t-1}+\delta_{t-1},\quad t\in \Z,$$
have a unique solution given by $U_t=W_t^{(1)}$ and the latter series is almost surely convergent. See \citep{Brandt} for a proof.
As a consequence, for any $t\in \Z$, $W_t^{(1)}<\infty$ a.s.

Finally, let us show that $W_t^{(4)}$ is finite a.s. 
Set  
$$W_{2,t}^{(4)}:=\sum_{s=1}^{\infty}\kappa_{t+s-1}\cdots\kappa_{t-1}= \sum_{s=1}^{\infty}\exp\left[\sum_{k=0}^s\log \kappa_{t+k-1}\right].$$
From {\bf A2} and the ergodic theorem, $\lim_{s\rightarrow \infty}\sum_{k=0}^s\log \kappa_{t+k-1}=-\infty$ a.s. and then 
$W_{2,t}^{(4)}$ is a.s. finite. Since from {\bf A3}, $\phi$ is a polynomial function such that $\phi(h)=O(h)$ in a neighborhood of $h=0$, 
we also deduce that $W_t^{(4)}$ is finite a.s.

Now for $C>0$, let $p_C:=\P\left(W_1^{(1)}\leq C,W_1^{(4)}\leq C,W_1^{(2)}\leq 1-1/C,W_1^{(3)}\leq 1-1/C\right)$. Then,
\begin{eqnarray*}
\lim_{C\rightarrow \infty}p_C&=&\P\left(W_1^{(1)}<\infty,W_1^{(4)}<\infty,W_1^{(2)}<1,W_1^{(3)}<1\right)\\
&=&\P\left(W_1^{(2)}<1,W_1^{(3)}<1\right)>0.
\end{eqnarray*}
There then exists $C>1$ s.t. $p_C>0$ which leads to the result.$\square$
\bigskip

We now consider the successive random times $\tau_0:=0<\tau_1<\tau_2<\cdots$ such that $\tau_i-\tau_{i-1}>h$ a.s. and 
$$W_{\tau_i}^{(1)}\leq C, W_{\tau_i}^{(2)}\leq 1-1/C,W_{\tau_i}^{(3)}\leq 1-1/C,W_{\tau_i}^{(4)}\leq C\mbox{ a.s.}$$
From the ergodic properties stated in Lemma (\ref{control}), the number of such random times is almost surely infinite.
 
We set $Z_{\omega,i}=\lambda_{\tau_i(\omega)}$ and $Z'_{\omega,i}=\lambda'_{\tau_i(\omega)}$.
We also set 
$$W_{\omega,i}=\frac{V\left(Z_{\omega,i}\right)+V\left(Z'_{\omega,i}\right)}{2}.$$

\begin{lem}\label{first}
Under the probability measure $\overline{\P}_{\omega}$, the three processes $\left(Z_{\omega,i}\right)_{i\geq 0}$, $\left(Z'_{\omega,i}\right)_{i\geq 0}$ $\left(\left(Z_{\omega,i},Z'_{\omega,i}\right)\right)_{i\geq 0}$ are time-inhomogeneous Markov chains. Moreover, for almost every $\omega$ and for $i\geq 2$,
$$\overline{\E}_{\omega}\left[V\left(Z_{\omega,i}\right)\vert Z_{\omega,i-1}\right]\leq (1-1/C)V\left(Z_{\omega,i-1}\right)+C,$$
$$\overline{\E}_{\omega}\left[V\left(Z'_{\omega,i}\right)\vert Z'_{\omega,i-1}\right]\leq (1-1/C)V\left(Z_{\omega,i-1}\right)+C$$
and 
$$\overline{\E}_{\omega}\left[V\left(Z_{\omega,1} \right)\right]\leq (1-1/C)V(s) +C,\quad \overline{\E}_{\omega}\left[V\left(Z'_{\omega,1} \right)\right]\leq (1-1/C)V(s) +C.$$
\end{lem}

\paragraph{Proof of Lemma \ref{first}}
The three processes $\left(Z_{\omega,i}\right)_{i\geq 0}$, $\left(Z'_{\omega,i}\right)_{i\geq 0}$ and $\left(Z'_{\omega,i}\right)_{i\geq 0}$ are subsequences of the Markov chains $\left(\lambda_t\right)_{t\geq 0}$, $\left(\lambda'_t\right)_{t\geq 0}$ and respectively 
$\left(\left(\lambda_t,\lambda'_t\right)\right)_{t\in \Z}$, they then keep the Markov property.
For proving the second part of the lemma, we use Lemma \ref{control} and the definition of the random times $\tau_i$, $i\geq 0$.
For $i\geq 2$, we have from {\bf A2},
\begin{eqnarray*}
\overline{\E}_{\omega}\left[V\left(Z_{\omega,i}\right)\vert Z_{\omega,i-1}\right]&\leq& \prod_{s=\tau_{i-1}(\omega)}^{\tau_i(\omega)-1}\gamma_s(\omega)V\left(Z_{\omega,i-1}\right)+W^{(1)}_{\tau_i(\omega)}\\
&\leq& W_{\tau_i(\omega)}^{(2)}V\left(Z_{\omega,i-1}\right)+W^{(1)}_{\tau_i(\omega)}\\
&\leq& (1-1/C)V\left(Z_{\omega,i-1}\right)+C.
\end{eqnarray*}
The same bound holds true for the other quantities.$\square$
\bigskip

Next we set
$$\rho_{\omega,1}=\inf\left\{i\geq 0\colon\; V\left(Z_{\omega,i}\right) +V\left(Z'_{\omega,i}\right) \leq C_1\right\}$$
and for $j\geq 2$ and a positive integer $d$,
$$\rho_{\omega,j}=\inf\left\{i\geq \rho_{\omega,j-1}+d\colon\; V\left(Z_{\omega,i}\right)+V\left(Z'_{\omega,i}\right) \leq C_1\right\}.$$
Let $W_{\omega,i}=\frac{V\left(Z_{\omega,i}\right)+V\left(Z'_{\omega,i}\right)}{2}$.
We also have 
$$\overline{\E}_{\omega}\left[W_{\omega,i}\vert Z_{\omega,i-1},Z'_{\omega,i-1}\right]\leq \kappa W_{\omega,i-1}+C.$$ 
Then,
$$\overline{\E}_{\omega}\left[W_{\omega,i}\vert W_{\omega,i-1}\right]\leq \kappa W_{\omega,i-1}+C.$$ 

We also set $\mathcal{F}_t=\sigma\left(\left(Y_t,Y'_t\right): t\geq 0\right)$ and $T_{\omega,j}=\tau_{\rho_{\omega,j}}(\omega)$ for $j\geq 1$. Observe that $T_{\omega,j}-1$ is a $\left(\mathcal{F}_t\right)_{t\geq 0}$-stopping time.
Moreover setting for $i\geq 1$, $\mathcal{G}_i^{\omega}=\mathcal{F}_{\tau_i(\omega)-1}$, the random variables $\rho_{\omega,j}$ are $\left(\mathcal{G}_i^{\omega}\right)_{i\geq 1}$-stopping times. Note also that $\mathcal{G}^{\omega}_{\rho_{\omega,j}}=\mathcal{F}_{T_{\omega,j}-1}$.

In what follows, we set $x=\frac{V(s)+V(s')}{2}$. In the following lemma, we obtain a control of some exponential moments for the stopping times $\rho_{\omega,j}$, $j\geq 1$.

\begin{lem}\label{annex}
Let $C_1=C(2 C+2)$, $\eta=\frac{2}{2-1/C}$ and $C_3=1+C^2+(C-1)C_1$. 
\begin{enumerate}
\item 
If $x>C_1$, we have $\overline{\E}_{\omega}\left[\eta^{\rho_{\omega,1}}\right]\leq x$.
\item
For all $i\geq 1$, $\overline{\E}_{\omega}\left[\eta^{\rho_{\omega,i+1}-\rho_{\omega,i}}\vert \mathcal{F}_{T_{\omega,i}-1}\right]\leq C_3\eta^d$.
\end{enumerate}
\end{lem}

\paragraph{Proof of Lemma \ref{annex}}

From Lemma \ref{first}, $\left(\left(Z_{\omega,i},Z'_{\omega_i}\right)\right)_{i\geq 0}$ is a Markov chain satisfying  drift conditions 
with fixed parameters (i.e. not depending on the index $i$). Lemma $3$ in \citep{DN} gives a similar control of such exponential moments for time-homogeneous Markov chains. However, their proof only uses the same kind of drift condition and it is also valid for time-inhomogeneous Markov chains, provided that the drift condition does not depends on the time index. Since their arguments are exactly the same, we omit the details.$\square$

\subsection{Lower bound for the probability (\ref{aborne})}\label{second}
We set 
$$M_n(\omega)=\sup\left\{i\geq 0\colon\; \tau_i(\omega)\leq n\right\}$$
From the ergodicity of the process $\left(U_t\right)_{t\in\Z}$ defined by $U_t=\left(W_t^{(i)}\right)_{1\leq i\leq 4}$, it is clear that 
$M_n(\omega)\rightarrow \infty$. 
If $d=d_n(\omega)\rightarrow \infty$ is a sequence such that $d_n(\omega)=o\left(M_n(\omega)\right)$,
we introduce the stopping time
$$\rho_{\omega}^{(n)}=\inf\left\{ i\geq 0\colon\; \left\vert \lambda_{\tau_i(\omega)}-\lambda'_{\tau_i(\omega)}\right\vert\leq C_1^{1/\alpha}\left(1-1/C\right)^{d_n}\right\}.$$
Our aim is to bound $\overline{\P}_{\omega}\left(\rho_{\omega}^{(n)}\geq M_n(\omega)\right)$. 

Consider a sequence $N_n(\omega)\rightarrow \infty$ and such that $N_n(\omega)=o\left(M_n(\omega)\right)$.
For simplicity of notations, we simply note $M_n(\omega),d_n(\omega)$ and $N_n(\omega)$ by $M_n,d_n$ and $N_n$.
We also note $\tau_{\rho_{\omega,j}+s}(\omega)$ by $T_{\omega,j,s}$. Note that 
$T_{\omega,j,0}=T_{\omega,j}$ which has been previously defined. 
We have the bound 
$$\overline{\P}_{\omega}\left(\rho_{\omega}^{(n)}\geq M_n\right)\leq \overline{\P}_{\omega}\left(\rho_{\omega,N_n+d_n}\geq M_n\right)+\overline{\P}_{\omega}\left(A_1^c\cap\cdots\cap A_{N_n}^c\right),$$
where for $1\leq j\leq N_n$, 
$$A_j=\left\{Y_t=Y'_t: T_{\omega,j}\leq t< T_{\omega,j,d_n}\right\}.$$
Indeed, on the event $A_j\cap\left\{\rho_{\omega}^{(n)}\geq M_n, \rho_{\omega,N_n+d_n}<M_n\right\}$,  we have from {\bf A1} and the definition of the random times $\tau_i$,
\begin{eqnarray*}
\left\vert \lambda_{T_{\omega,j,d_n}}-\lambda'_{T_{\omega,j,d_n}}\right\vert&\leq &\prod_{T_{\omega,j}\leq t<T_{\omega,j,d_n}}\kappa_t(\omega)\left\vert \lambda_{T_{\omega,j}}-\lambda'_{T_{\omega,j}}\right\vert\\
&\leq& C_1^{1/\alpha}\prod_{1\leq s\leq d_n}\prod_{T_{\omega,j,s-1}\leq t<T_{\omega,j,s}}\kappa_t(\omega)\\
&\leq& C_1^{1/\alpha}\prod_{1\leq s\leq d_n}W^{(3)}_{T_{\omega,j,s}}(\omega)\\
&\leq & C_1^{1/\alpha} \left(1-1/C\right)^{d_n}
\end{eqnarray*}
and the intersection of $A_j$ with the event $\left\{\rho_{\omega}^{(n)}\geq M_n\right\}$ is then empty. In the previous bounds, we used the inequality
$$\left\vert \lambda_{T_{\omega,j}}-\lambda'_{T_{\omega,j}}\right\vert\leq \left(V\left(\lambda_{T_{\omega,j}}\right)+V\left(\lambda'_{T_{\omega,j}}\right)\right)^{1/\alpha}\leq C_1^{1/\alpha}.$$
Next, from {\bf A3}, the function $\phi$ is non-decreasing and if $\overline{C}_1=\max\left(1,C_1^{1/\alpha}\right)^{d(\phi)}$ where 
$d(\phi)$ denotes the degree of the polynomial function $\phi$, we have the bound $\phi\left(C_1^{1/\alpha}h\right)\leq \overline{C}_1\phi(h)$ for all $h\geq 0$.
We then get from {\bf A3}, 
\begin{eqnarray*}
\overline{\P}_{\omega}\left(A_j\vert \mathcal{F}_{T_{\omega,j}-1}\right)&\geq & \overline{\P}_{\omega}\left(Y_{T_{\omega,j}+s}=Y'_{T_{\omega,j}+s}; s\geq 0\vert \mathcal{F}_{T_{\omega,j}-1}\right)\\
&=& \lim_{\ell\rightarrow\infty}\prod_{s=0}^{\ell}\overline{\P}_{\omega}\left(Y_{T_{\omega,j}+s}=Y'_{T_{\omega,j}+s}\vert \mathcal{F}_{T_{\omega,j}-1}, Y_{T_{\omega,j}+i}=Y'_{T_{\omega,j}+i}; 0\leq i\leq s-1\right)\\
&\geq& \lim_{\ell\rightarrow \infty}\prod_{s=0}^{\ell}\exp\left(-\phi\left(\prod_{i=0}^{s-1}\kappa_{T_{\omega,j}+i}(\omega)\left\vert \lambda_{T_{\omega,j}}-\lambda'_{T_{\omega,j}}\right\vert\right)\right)\\
&\geq& \exp\left(-\phi\left(C_1^{1/\alpha}\right)-\overline{C}_1W^{(4)}_{T_{\omega,j}}(\omega)\right)\\
&\geq& \exp\left(-\phi\left(C_1^{1/\alpha}\right)-\overline{C}_1C\right):=C_2.
\end{eqnarray*}
We deduce the bound 
$$\overline{\P}_{\omega}\left(A_1^c\cap\cdots\cap A_{N_n}^c\right)\leq \left(1-C_2\right)^{N_n}.$$
Using Lemma \ref{annex}, we also get 
\begin{eqnarray*}
\overline{\P}_{\omega}\left(\rho_{\omega,N_n+d_n}\geq M_n\right)&\leq& \eta^{-M_n}\overline{\E}_{\omega}\left[\eta^{\rho_{\omega,1}+\sum_{i=2}^{N_n+d_n}\left(\rho_{\omega,i}-\rho_{\omega,i-1}\right)}\right]\\
&\leq& \eta^{-M_n}\overline{\E}_{\omega}\left[\eta^{\rho_{\omega,1}}\right]C_3^{N_n+d_n-1}\eta^{d_n(d_n+N_n-1)}.
\end{eqnarray*}
Note that
$$\overline{\E}_{\omega}\left[\eta^{\rho_{\omega,1}}\right]\leq 1+\overline{\E}_{\omega}\left[\eta^{\rho_{\omega,1}}\mathds{1}_{x>C_1}\right]\leq 1+x.$$
We then obtain the following bound.
\begin{equation}\label{bound1}
\overline{\P}_{\omega}\left(\rho_{\omega}^{(n)}\geq M_n\right)\leq \left(1-C_2\right)^{N_n}+(1+x)C_3^{d_n+N_n-1}\eta^{d_n(N_n+d_n-1)-M_n}.
\end{equation}
Next, writting $\phi(h)=\sum_{j=1}^{d(\phi)}\phi_jh^j$, we set $C_4=\overline{C}_1\phi(1)+C\overline{C}_1$.
Setting $S_{\omega}^{(n)}=\tau_{\rho_{\omega}^{(n)}}(\omega)$, on the event $\left\{\rho_{\omega}^{(n)}<M_n\right\}=\left\{S_{\omega}^{(n)}<\tau_{M_n(\omega)}(\omega)\right\}$, we have 
\begin{eqnarray*}
&&\overline{\P}_{\omega}\left(Y_{S_{\omega}^{(n)}+i}=Y'_{S_{\omega}^{(n)}+i}; i\geq 0 \vert \mathcal{F}_{S_{\omega}^{(n)}}\right)\\&=& \prod_{i=0}^{\infty}\overline{\P}_{\omega}\left(Y_{S_{\omega}^{(n)}+i}=Y'_{S_{\omega}^{(n)}+i}\vert \mathcal{F}_{S_{\omega}^{(n)}-1},Y_{S_{\omega}^{(n)}+j}=Y'_{S_{\omega}^{(n)}+j}; 0\leq j\leq i-1\right)\\
&\geq& \prod_{i=0}^{\infty}\exp\left(-\phi\left(\kappa_{S_{\omega}^{(n)}+i-1}(\omega)\cdots \kappa_{S_{\omega}^{(n)}}(\omega)\left\vert \lambda_{S_{\omega}^{(n)}}-\lambda'_{S_{\omega}^{(n)}}\right\vert\right)\right)\\
&\geq& \exp\left(-\phi\left(C_1^{1/\alpha}(1-1/C)^{d_n}\right)-\overline{C}_1(1-1/C)^{d_n}W^{(4)}_{S_{\omega}^{(n)}}\right)\\
&\geq & \exp\left(-C_4(1-1/C)^{d_n}\right).
\end{eqnarray*}
Moreover, using the inequality $h\leq \phi_1^{-1}\phi(h)$, we have on the event $\left\{Y_{S_{\omega}^{(n)}+i}=Y'_{S_{\omega}^{(n)}+i}; i\geq 0\right\}$,
\begin{eqnarray*}
\sum_{i=1}^{\infty}\left\vert \lambda_{S_{\omega}^{(n)}+i}-\lambda'_{S_{\omega}^{(n)}+i}\right\vert &\leq& \sum_{i=1}^{\infty}\prod_{\ell=0}^{i-1}\kappa_{S_{\omega}^{(n)}+\ell}(\omega) \left\vert \lambda_{S_{\omega}^{(n)}}-\lambda'_{S_{\omega}^{(n)}}\right\vert\\
&\leq&C_1^{1/\alpha}\phi_1^{-1}(1-1/C)^{d_n}W^{(4)}_{S_{\omega}^{(n)}}\\
& \leq& CC_1^{1/\alpha}\phi_1^{-1}(1-1/C)^{d_n}.
\end{eqnarray*}
We then deduce the following result.

\begin{prop}\label{inter}
There exist $\widetilde{C}>0$ and $\overline{\widetilde{\rho}}\in (0,1)$, only depending on $C,C_1,\phi$ and $\alpha$ such that
$$\overline{\P}_{\omega}\left(Y_{n+i}=Y'_{n+i}; i\geq 0 \mbox{ and } \sum_{i=1}^{\infty}\left\vert \lambda_{n+i}-\lambda'_{n+i}\right\vert \leq \widetilde{C}\widetilde{\rho}^{\sqrt{M_n(\omega)}}\right)\geq 1-\widetilde{C}(1+x)\widetilde{\rho}^{\sqrt{M_n(\omega)}}.$$
\end{prop}

\paragraph{Proof of Proposition \ref{inter}}
We choose $d_n=[\sqrt{M_n(\omega)}/2]$ and $N_n=[\sqrt{M_n(\omega)}/2]$ where $[x]$ denotes is the integer part of a real number $x$.
There exist $\widetilde{C}>0$ and $\widetilde{\rho}\in (0,1)$ large enough for getting the bounds $CC_1^{1/\alpha}\phi_1^{-1}(1-1/C)^{d_n}\leq \widetilde{C}\widetilde{\rho}^{\sqrt{M_n(\omega)}}$ and
$$\exp\left(-C_4(1-1/C)^{d_n}\right)\cdot\left(1-\left(1-C_2\right)^{N_n}-(1+x)C_3^{d_n+N_n-1}\eta^{d_n(N_n+d_n-1)-M_n(\omega)}\right)\geq 1-\widetilde{C}(1+x)\widetilde{\rho}^{\sqrt{M_n(\omega)}}.$$
Setting 
$$A_n=\left\{Y_{n+i}=Y'_{n+i}; i\geq 0 \mbox{ and } \sum_{i=1}^{\infty}\left\vert \lambda_{n+i}-\lambda'_{n+i}\right\vert \leq CC_1^{1/\alpha}\phi_1^{-1}(1-1/C)^{d_n}\right\},$$
$$B_n=\left\{Y_{S_{\omega}^{(n)}+i}=Y'_{S_{\omega}^{(n)}+i}; i\geq 0 \mbox{ and } \sum_{i=1}^{\infty}\left\vert \lambda_{S_{\omega}^{(n)}+i}-\lambda'_{S_{\omega}^{(n)}+i}\right\vert \leq CC_1^{1/\alpha}\phi_1^{-1}(1-1/C)^{d_n}\right\},$$
we have, using our previous computations, (\ref{bound1}) and the definition of $\widetilde{C}$ and $\widetilde{\rho}$,
\begin{eqnarray*}
\overline{\P}_{\omega}(A_n)&=&\overline{\P}_{\omega}\left(A_n\cap\left\{S_{\omega}^{(n)}< \tau_{M_n(\omega)}(\omega)\right\}\right)\\
&\geq& \overline{\P}_{\omega}\left(B_n\cap\left\{S_{\omega}^{(n)}< \tau_{M_n(\omega)}(\omega)\right\}\right)\\
&\geq& \overline{\E}_{\omega}\left[\overline{\P}_{\omega}\left(B_n\vert \mathcal{F}_{S_{\omega}^{(n)}-1}\right)\mathds{1}_{\left\{S_{\omega}^{(n)}< \tau_{M_n(\omega)}(\omega)\right\}}\right]\\
&\geq& \exp\left(-C_4(1-1/C)^{d_n}\right)\overline{\P}_{\omega}\left(S_{\omega}^{(n)}< \tau_{M_n(\omega)}(\omega)\right)\\
&\geq& 1-\widetilde{C}(1+x)\widetilde{\rho}^{\sqrt{M_n(\omega)}}.
\end{eqnarray*}
This leads to the proposed lower bound.$\square$

\subsection{Convergence of the backward iterations and proof of Theorem \ref{main}}\label{third}
Next, we consider on $F$ the metric $\Delta(s,s')=\vert s-s'\vert\wedge 1$ and the corresponding Wasserstein metric of order $1$,
$$\mathcal{W}_1\left(\mu,\nu\right)=\inf\left\{\int \Delta(s,s')\gamma(ds,ds')\right\},$$
where the infimum is on the set of probability measures $\gamma$ on $F\times F$ possessing marginals $\mu$ and $\nu$.
In what follows, we denote by $\delta_s$ the Dirac mass at point $s\in F$. We also remind that if $\mu$ is a probability measure on $F$ and $P$ is a Markov kernel on $F$, the probability $\mu P$ is defined by $\mu P(A)=\int \mu(ds)P(s,A)$ for any $A\in \mathcal{B}(F)$.

\begin{prop}\label{inter2}
There exist $\overline{C}>0$ and $\overline{\rho}\in (0,1)$, only depending on $C,C_1,\phi$ such that
$$\mathcal{W}_1\left(\delta_sP_{X_0(\omega)}\cdots P_{X_{n-1}(\omega)},\delta_{s'}P_{X_0(\omega)}\cdots P_{X_{n-1}(\omega)}\right)\leq \overline{C}\left(1+V(s)+V(s')\vert\right)\rho^{\sqrt{M_n(\omega)}}.$$
\end{prop}

\paragraph{Proof of Proposition \ref{inter2}}
We use Proposition \ref{inter}. Since 
$$\Delta(s,s')\leq \widetilde{C}\widetilde{\rho}^{\sqrt{M_n(\omega)}}+\mathds{1}_{\left\{\vert s-s'\vert>\widetilde{C}\widetilde{\rho}^{\sqrt{M_n(\omega)}}\right\}},$$
we have 
\begin{eqnarray*}
\mathcal{W}_1\left(\delta_sP_{X_0(\omega)}\cdots P_{X_{n-1}(\omega)},\delta_{s'}P_{X_0(\omega)}\cdots P_{X_{n-1}(\omega)}\right)&\leq& 
\overline{\E}_{\omega}\left[\Delta\left(\lambda_n,\lambda'_n\right)\right]\\
&\leq&\widetilde{C}\widetilde{\rho}^{\sqrt{M_n(\omega)}}+\overline{\P}_{\omega}\left(\left\vert\lambda_n-\lambda'_n\right\vert>\widetilde{C}\widetilde{\rho}^{\sqrt{M_n(\omega)}}\right)\\
&\leq& \widetilde{C}\widetilde{\rho}^{\sqrt{M_n(\omega)}}+\widetilde{C}(1+x)\widetilde{\rho}^{\sqrt{M_n(\omega)}}.
\end{eqnarray*}
We then get the result by setting $\overline{\rho}=\widetilde{\rho}$ and $\overline{C}=2\widetilde{C}$.$\square$
\bigskip

\begin{prop}\label{inter3}
Let Assumptions {\bf A1-A3} hold true. There then exists a unique process $\left(\pi_t\right)_{t\in\Z}$ of identically distributed random probability measures on $F$ such that and such that $\pi_t P_{X_t}=\pi_{t+1}$ a.s. 
Moreover, almost surely, for any $s\in F$, 
$$\lim_{n\rightarrow \infty}\mathcal{W}_1\left(\delta_sP_{X_{t-n}}\cdots P_{X_{t-1}},\pi_t\right)=0.$$
\end{prop}

\paragraph{Proof of Proposition \ref{inter3}}
For simplicity, we now work on the canonical space $\Omega=L^{\Z}$ and we assume that $X_t(\omega)=\omega_t$. We will simply denote 
$P_{X_t(\omega)}$ by $P^{\omega_t}$. We are going to show that there exists a random probability measure $\mu^{\omega}$ such that 
for any $s\in F$,
$$\lim_{n\rightarrow \infty}\mathcal{W}_1\left(\delta_sP^{\omega_{-n}}\cdots P^{\omega_{-1}},\mu^{\omega}\right)=0\mbox{ a.s.}.$$
First, it is easily seen from Proposition \ref{inter2} that for $s,s'\in F$,  
\begin{equation}\label{target}
\mathcal{W}_1\left(\delta_sP^{\omega_{-n}}\cdots P^{\omega_{-1}},\delta_{s'}P^{\omega_{-n}}\cdots P^{\omega_{-1}}\right)\leq \overline{C}\left(1+V(s)+V(s')\vert\right)\overline{\rho}^{\sqrt{M_n(\theta^{-n}\omega)}},
\end{equation}
where $\theta\colon\;L^{\Z}\rightarrow L^{\Z}$ is the shift operator defined by $\theta\omega=\left(\omega_{t+1}\right)_{t\in\Z}$.
We set $U_t=\left(W_t^{(i)}\right)_{1\leq i\leq 4}$ and $\mathcal{C}=[0,C]\times [0,1-1/C]^2\times [0,C]$ where all the quantities are defined in Lemma \ref{control}. Set also $L_n(\omega)=M_n\left(\theta^{-n}\omega\right)$ and . Note that $\tau_i\left(\theta^{-n}\omega\right)-n$ are precisely 
the successive time points $t$ between $-n$ and $0$, distant at least of $h$ units of times and such that $U_t(\omega)\in\mathcal{C}$.  Then, from the ergodicity of the process $\left(U_t\right)_{t\in\Z}$, $\lim_{n\rightarrow \infty} L_n(\omega)=\lim_{n\rightarrow \infty} M_n(\omega)=\infty$.
Note that from (\ref{target}), if the random probability measure $\mu^{\omega}$ exists, it cannot depend on $s$.
Let $t_n(\omega):=\tau_1\left(\theta^{-n}\omega\right)-n$ be now the first time point $t\geq -n$ such that $U_t(\omega)\in\mathcal{C}$.
Now for $m\geq n$, and $s_1,s_2\in F$,

\begin{eqnarray*}
&&\mathcal{W}_1\left(\delta_{s_1} P^{\omega_{-m}}\cdots P^{\omega_{-1}},\delta_{s_2}P^{\omega_{-n}}\cdots P^{\omega_{-1}}\right)\\&\leq &\int\int \delta_{s_1}P^{\omega_{-m}}\cdots P^{\omega_{t_n(\omega)-1}}(ds_3)\delta_{s_2}P^{\omega_{-n}}\cdots P^{\omega_{t_n(\omega)-1}}(ds_4)W_1\left(\delta_{s_3}P^{\omega_{t_n(\omega)}}\cdots P^{\omega_{-1}},\delta_{s_4}P^{\omega_{t_n(\omega)}}\cdots P^{\omega_{-1}}\right)\\
&\leq&\int\int \delta_{s_1} P^{\omega_{-m}}\cdots P^{\omega_{t_n(\omega)}}(ds_3)\delta_{s_2}P ^{\omega_{-n}}\cdots P^{\omega_{t_n(\omega)}}(ds_4)\overline{C}\left(1+V(s_3)+V(s_4)\right)\rho^{\sqrt{L_n(\omega)-1}}.
\end{eqnarray*}
In the last bound, we used (\ref{target}) and the inequality $L_{-t_n(\omega)}(\omega)=L_n(\omega)-1$. 
Note the from the drift condition {\bf A2} and the definition of $t_n(\omega)$, we have
$$\delta_{s_1} P^{\omega_{-m}}\cdots P^{\omega_{t_n(\omega)}}V\leq (1-1/C)V(s_1)\vert+C\leq V(s_1)+C.$$
We then obtain
$$\mathcal{W}_1\left(\delta_{s_1} P^{\omega_{-m}}\cdots P^{\omega_{-1}},\delta_{s_2}P^{\omega_{-n}}\cdots P^{\omega_{-1}}\right)
\leq \overline{C}\left(1+2C+V(s_1)+V(s_2)\right)\rho^{\sqrt{L_n(\omega)-1}}.$$
Setting $s_1=s_2=s$, one can see that the sequence $\left(\mu_{n,s}^{\omega}\right)_{n\geq 1}$ defined by 
$\mu_{n,s}^{\omega}=\delta_sP^{\omega_{-n}}\cdots P^{\omega_{-1}}$  is a Cauchy sequence in the complete space of probability measures 
endowed with the metric $\mathcal{W}_1$. 
We can then define a limit $\mu^{\omega}=\lim_{n\rightarrow \infty}\delta_s P^{\omega_{-n}}\cdots P^{\omega_{-1}}$. 
As previously mentioned, this probability measure does not depend on $s$.
We then set $\pi_t(\omega)=\mu^{\theta^t\omega}$. Clearly, the sequence $\left(\pi_t\right)_{t\in\Z}$ is a stationary sequence 
of random probability measures.
It is only necessary to check the equality $\mu^{\theta^t\omega} P^{\omega_t}=\mu^{\theta^{t+1}\omega}$ for $t=0$. To this end, let $g:F\rightarrow \R$ be a bounded and Lipschitz function. Then $g$ is also a Lipschitz function from $\left(L,d\right)$ to $\R$. Moreover, $\P^{\omega_1}g$ is also a Lipschitz function from $\left(L,d\right)$ to $\R$. Indeed, it is a bounded function and from {\bf A2-A3}, we have
\begin{eqnarray*}
\left\vert P^{\omega_1}g(s)-P^{\omega_1}g(s')\right\vert&\leq& L(g)\kappa\vert s-s'\vert+\Vert g\Vert_{\infty}d_{TV}\left(p(\cdot\vert s),p(\cdot\vert s')\right)\\
&\leq& \left(L(g)\kappa+K\Vert g\Vert_{\infty}\right)\vert s-s'\vert,
\end{eqnarray*}
where $L(g)$ denotes the Lipschitz constant of $g$ and $\Vert g\Vert_{\infty}=\sup_{s\in F}\left\vert g(s)\right\vert$.
Since convergence in Wasserstein metric $\mathcal{W}_1$ entails convergence of the integrals of Lipschitz functions, we get  
$$\mu^{\omega} P^{\omega_0} g=\lim_{n\rightarrow \infty}\delta_z P^{\omega_{-n}}\cdots P^{\omega_{-1}}(P^{\omega_0}g)=\lim_{n\rightarrow \infty}\delta_z P^{\omega_{-n}}\cdots P^{\omega_0}g=\mu^{\theta\omega}g\mbox{ a.s.}$$
We then deduce the equality $\mu^{\omega}P^{\omega_0}=\mu^{\theta \omega}$.

Next, we show uniqueness. Let $(\overline{\pi}_t)_{t\in\Z}$ be another process of identically distributed random variables and such that $\overline{\pi}_tP_{X_t}=\overline{\pi}_{t+1}$.
We have $\overline{\pi}_0=\overline{\pi}_{-n}P_{X_{-n}}\cdots P_{X_{-1}}$ a.s. 
and for a given $s\in F$, 
\begin{eqnarray*}
\E\vert \overline{\pi}_0 g -\delta_z P_{X_{-n}}\cdots P_{X_{-1}} g\vert&\leq &\E\int_{s'} \overline{\pi}_{-n}(ds')\vert
\delta_{s'}P_{X_{-n}}\cdots P_{X_{-1}} g-\delta_s P_{X_{-n}}\cdots P_{X_{-1}} g\vert\\
&\leq & 2\E\overline{\pi}_0(\{\vert s'\vert>M\})\\
& & {} + L(g)\E\sup_{\vert s'\vert\leq M}\mathcal{W}_1\left(\delta_{s'} P_{X_{-n}}\cdots P_{X_{-1}}-\delta_s P_{X_{-n}}\cdots P_{X_{-1}}\right).
\end{eqnarray*}
The first term of this last bound can be made arbitrarily small when $M$ is large, using Lebesgue's Theorem. Moreover, for a fixed positive $M$, the second term goes to $0$ as $n$ goes to infinity, using 
(\ref{target}) and Lebesgue's theorem. Since, $\delta_s P_{X_{-n}}\cdots P_{X_{-1}} g\rightarrow \pi_0 g$, we deduce that $\overline{\pi}_0 g=\pi_0 g$ a.s. and since $g$ is an arbitrary Lipschitz function, $\overline{\pi}_0=\pi_0$ a.s.$\square$

\subsection{Proof of Theorem \ref{main}}
If $R^{\omega_0}=R_{X_0(\omega)}$ be the Markov kernel defined just after (\ref{dynamic}). As explained, conditional on $X$, any solution $(Y_t,\lambda_t)$ of (\ref{dynamic}) is a time-inhomogeneous Markov chain with transition kernels $\left(R^{\omega_t}\right)$. It is easily seen that if $\left(\nu_t\right)_{t\in\Z}$ is a sequence of identically distributed random measures on $E\times F$ and such that $\nu_t(\omega)R^{\omega_t}=\nu_{t+1}(\omega)$ a.s. then
\begin{equation}\label{details1}
\nu_t(\omega)(dy,ds)=p(dy\vert s)\mu^{\theta^t\omega}(ds)\mbox{ a.s.}$$
\end{equation}
Indeed, if for $B\in\mathcal{B}(F)$, $\nu_t^{(2)}(\omega)(B)=\nu_t(\omega)\left(E\times B\right)$ (the second marginal of $\nu_t(\omega)$), the invariance relation leads to 
$$\nu_t(\omega)(dy',ds')=p(dy'\vert s')\nu_t^{(2)}(ds'),\quad \nu_t^{(2)}(\omega)P^{\omega_t}(ds')=\nu^{(2)}_{t+1}(\omega)(ds').$$
From, the uniqueness property in Proposition \ref{inter3}, we get (\ref{details1}). In what follows, set $\nu^{\omega}=\nu_0(\omega)$.
Let us now prove Theorem \ref{main}. For the existence part, we consider the finite-dimensional distributions 
$$\zeta_{u,t}^{\omega}(dy_u,d_{s_u},\cdots,dy_t,ds_t)=\nu^{\theta^{u}\omega}\left(dy_{u},ds_{u}\right)\prod_{i=u}^{t-1}R^{\omega_i}\left((y_i,s_i),(dy_{i+1},ds_{i+1})\right),$$
for $u\leq t$ in $\Z$. Using Kolmogorov's extension theorem, there exists a unique probability measure $\zeta^{\omega}$ on $(E\times F)^{\Z}$ compatible with such a family.
On $\widetilde{\Omega}=(E\times F)^{\Z}\times \Omega$, the probability measure $d\mathbb{Q}\left((y,s),\omega\right)=d\zeta^{\omega}\left((y,s)\right)d\P(\omega)$ is solution of (\ref{dynamic}), in the sense that if $Y_t((y,s),\omega)=y_t$ and $\lambda_t\left((y,s),\omega\right)=s_t$, we have
$$\mathbb{Q}\left(Y_t\in A\vert X,Y_{t-1},\lambda_t,Y_{t-2},\lambda_{t-1}\ldots\right)=p(A\vert \lambda_t),\quad\lambda_t=f\left(\lambda_{t-1},Y_{t-1},X_{t-1}\right).$$
Stationarity of such a solution results from the equalities
$$\zeta_{u,t}^{\theta^n \omega}=\zeta_{u+n,t+n}^{\omega},$$
for $u\leq t$ in $\Z$ and $n\in \N$. Let us now show uniqueness. If $\left((Y_t,\lambda_t,X_t)\right)_{t\in\Z}$ is a stochastic process satisfying (\ref{dynamic}), we will have almost surely,
$$\P\left(Y_t\in A,\lambda_t\in B\vert X\right)=\nu^{\theta^t\omega}(A\times B),\quad (A,B)\in \mathcal{B}(E)\times \mathcal{B}(F).$$
Indeed, the (conditional) marginal distribution is a random probability measure $\overline{\nu}_t$ such that $\overline{\nu}_t(\omega)R^{\theta^t\omega}=\overline{\nu}_{t+1}(\omega)$ for almost every $\omega$.
Hence we have $\overline{\nu}_t(\omega)=\nu_t(\omega)$ for almost every $\omega$. The probability distribution of $\left(\left(Y_t,\lambda_t\right)_{t\in\Z},\left(X_t\right)_{t\in\Z}\right)$ then coincides with $\mathbb{Q}$.

We know prove ergodicity of the unique stationary solution.
The uniqueness property derived previously can be used for proving ergodicity. 
We consider a Markov kernel $Q$ on $(E\times F)\times \Omega$ defined by
$$Q\left((z,\omega),(dz',d\omega')\right)=R^{\omega}(z,dz')\delta_{\theta\omega}(d\omega').$$
Note that the measure $\gamma((dy,ds),d\omega)=p(dy\vert s)\mu^{\omega}(ds)\P(d\omega)$ 
is invariant for $Q$. A Markov chain with transition $Q$ is usually referred as a skew-product Markov chain. See \citep{Kifer} or \citep{Orey}.
Moreover, we are going to show that $\gamma$ is the unique $Q-$invariant probability measure for which the second marginal, denoted by $\gamma_2$, is absolutely continuous with respect to $\P$, the distribution of the environment.
From the Radon-Nikodym theorem, there exists a measurable function $f\colon\;\Omega\rightarrow \R_+$ such that $d\gamma_2(\omega)=f(\omega)d\P(\omega)$.
For such a measure $\gamma$, we use a measure disintegration $\gamma(dz,d\omega)=\gamma_c(\omega,dz)f(\omega)\P(d\omega)$ where $\gamma_c$ is a probability kernel from $\Omega$ to 
$E\times F$. See for instance \citep{Kal} for a proof of existence for such probability kernel.
Necessarily, we get from the invariance equation, for $\P-$almost $\omega$,
$$\int \gamma_c(\omega,dz)R^{\omega}(z,dz')=\gamma_c\left(\theta\omega,dz'\right)\mbox{ and }f(\theta \omega)=f(\omega).$$
Ergodicity of $\P$ entails that $f=1$ $\P-$a.s. Moreover $\gamma_c\left(\theta^t,\cdot\right)$ has a distribution not depending on $t$ and then coincides $\P-$a.s. with $\nu_t$, using a uniqueness property proved just before. 
Suppose now that $\gamma$ is not ergodic. There then exists a measurable set $A$ such that $Q((z,\omega),A)=1$ for $\gamma-$almost $(z,\omega)\in A$ and $\gamma(A)\in (0,1)$.
From the invariance of $\gamma$, we also have $Q((z,\omega),A^c)=1$ for $\gamma-$almost $(z,\omega)\in A^c$. The measure $\gamma_A(B)=\gamma(A\cap B)/\gamma(A)$ is then another probability measure invariant for $Q$ and such that its second marginal is absolutely continuous w.r.t. $\P$. 
This contradicts the uniqueness property. Hence $\gamma$ is ergodic. Since $\gamma$ is the probability distribution of $\left(Y_0,\lambda_0,(X_t)_{t\in\Z}\right)$, we easily deduce the ergodicity of the process $\left((Y_t,\lambda_t,X_t)\right)_{t\in\Z}$.$\square$

\section{Examples}\label{S3}

We now give many time series models satisfying Assumptions {\bf A1-A3}. 
First we give a simple sufficient condition for checking the drift condition {\bf A2}.
For simplicity, we say that a function $\kappa\colon\;L\rightarrow \R$ is in $\mathcal{M}_{\log}$ if it is measurable and if $\E\log^{+}\left(\vert \kappa(X_0)\vert\right)<\infty$. In what follows, $i\in \{0,1\}$.

\begin{description}
\item[A2(i)] There exist functions $\kappa,\widetilde{\kappa},\widetilde{\delta}:F\rightarrow \R_+$ in $\mathcal{M}_{\log}$ such that 
$\E\log\kappa(X_0)<0$ and 
for $(s,y,x)\in F\times E\times L$,
$$\left\vert f(s,y,x)\right\vert\leq \kappa(x)\vert s\vert +\widetilde{\kappa}(x)\vert y\vert^i+\widetilde{\delta}(x).$$
\end{description}

The proof of the next result is straightforward.

\begin{prop}\label{pratique1}
Suppose that for some $i\in\{1,2\}$, Assumption {\bf A2(i)} holds true and there exists $D>0$ such that for any $s\in F$, $\int \vert y\vert^i p(dy\vert s)\leq \vert s\vert+D$. Condition {\bf A2} is then satisfied with $V(s)=1+\vert s\vert$. 
\end{prop}

We also provide a general result for checking {\bf A1-A2} when the latent process $\left(\lambda_t\right)_{t\in\Z}$ satisfies a threshold dynamic.

\begin{prop}\label{pratique2}
Suppose that the assumptions of Proposition \ref{pratique1} hold true and assume that $E\subset\R$. Assume furthermore that there exists functions $\kappa_j,\widetilde{\kappa}_j,\gamma_j \in \mathcal{M}_{\log}$
and some intervals $I(x)$, $x\in I$, of the real line such that
$$f(s,y,x)=\left\{\begin{array}{c}\kappa_1(x)s+\widetilde{\kappa}_1(x)y^i+\gamma_1(x),\mbox{ if } y\in I(x)\\
\kappa_2(x)s+\widetilde{\kappa}_2(x)y^i+\gamma_2(x),\mbox{ if } y\notin I(x)\end{array}\right.$$
\begin{enumerate}
\item
If $\E\log\left(\kappa(X_0)+\widetilde{\kappa}(X_0)\right)<0$ with $\kappa(x)=\max\left\{\vert \kappa_1(x)\vert, \vert \kappa_2(x)\vert\right\}$ and $\widetilde{\kappa}(x)=\max\left\{\vert \widetilde{\kappa}_1(x)\vert, \vert \widetilde{\kappa}_2(x)\vert\right\}$, conditions {\bf A1-A2} are fulfilled with $V(s)=1+\vert s\vert$.
\item
If for every $x\in L$, $I(x)$ is a bounded interval, the same conclusions hold true as soon as 
$$\E\log\left(\kappa(X_0)\right)<0\mbox{ and } \E\log\left(\left\vert \kappa_2(X_0)\right\vert+\left\vert \widetilde{\kappa}_2(X_0)\right\vert\right)<0.$$
\end{enumerate}
\end{prop}

\paragraph{Notes}

\begin{enumerate}
\item
When $\kappa_j$ and $\widetilde{\kappa}_j$ are deterministic for $j=1,2$, we simply have to assume the condition 
$$\max\{\vert \kappa_1\vert,\vert \kappa_2\vert\}+\max\left\{\vert \widetilde{\kappa}_1\vert,\vert \widetilde{\kappa}_2\vert\right\}<1$$
for the first point of Proposition \ref{pratique2} while for the second point, we only need the conditions
$$\vert \kappa_1\vert<1,\quad \vert \kappa_2\vert+\left\vert\widetilde{\kappa}_2\right\vert<1.$$
\item
Our framework allows random coefficients that depend on the exogenous covariates. Note that this point is important if we want to take in account of some interactions between lag values of the response and the covariates. 
For instance, $\widetilde{\kappa}_j(x)y^i=\widetilde{\kappa}_j xy^i$ for deterministic coefficients $\widetilde{\kappa}_j$, $j=1,2$.
\item
Note that from Jensen'inequality, any condition of type $\E\log\left(\kappa(X_0)\right)<0$ is satisfied as soon as $\E \kappa(X_0)<1$.
\end{enumerate}

For the various examples given in the rest of this section, {\bf A3} will be the main assumption to check. Assumptions {\bf A1-A2} can be checked using Proposition \ref{pratique1} or Proposition \ref{pratique2} for instance.

\subsection{Categorical time series}

We first consider binary processes with $E=\{0,1\}$ and $F=\R$.  For a cdf $F$ on $\R$, we set $p(1\vert s)=1-p(0\vert s)=F(s)$ and assume that
$$\P\left(Y_t\in A\vert X,Y_{t-1},\lambda_t,Y_{t-2},\lambda_{t-1}\ldots\right)=F\left(\lambda_t\right),\quad \lambda_t=f\left(\lambda_{t-1},Y_{t-1},X_{t-1}\right).$$
Here, we will enlighten that the semi-contractivity assumption {\bf A1} is sufficient to also get the drift condition {\bf A2}.
We state our result for the two main cases considered in the literature, when $F$ is the standard Gaussian c.d.f. (probit autoregressive model) or $F(s)=\frac{\exp(s)}{1+\exp(s)}$ (logistic autoregressive model).

\begin{prop}\label{binary}
Suppose that Assumption {\bf A1} holds true and that there exists $s_0\in F$ such that $x\mapsto f(s_0,y,x)$ is in $\mathcal{M}_{\log}$ for $y=0,1$. Then the conclusions of Theorem \ref{main} are valid 
for the probit or the logistic autoregressive model.
\end{prop}

\paragraph{Note.} Binary and categorical time series of this type have been introduced in the applied econometrics literature. See for instance \citep{kauppi}, \citep{russell} or \citep{rydberg}. In \citep{Truquet}, conditions ensuring existence of a stationary and ergodic solution when covariates are included in such dynamics are given. However, at least for the probit/logistic model and the multinomial autoregressions discussed below, we obtain here sharper results. Indeed, \citep{Truquet} used uniform Lipschitz type properties with respect to the covariate process $X$. As a consequence, it is not possible to consider the simple model 
$$\lambda_t=\kappa \lambda_{t-1}+\widetilde{\kappa}X_{t-1}Y_{t-1}+\gamma,$$
when the process $X$ is unbounded. In contrast, this case is covered by our result and the conditions $\vert \kappa\vert<1$ and 
$\E\log^{+}\vert X_0\vert<\infty$ are sufficient for ensuring existence and uniqueness of a stationary solution.

\paragraph{Proof of Proposition \ref{binary}}
We first check {\bf A2} for $V(s)=1+\vert s\vert$. This is automatic using {\bf A1} and the additional assumption since 
$$\left\vert f(s,y,x)\right\vert\leq \max_{y\in\{0,1\}}\left\vert f(s_0,y,x)\right\vert+\kappa(x)\vert s-s_0\vert\leq \kappa(x)\vert s\vert+\max_{y\in\{0,1\}}\left\vert f(s_0,y,x)\right\vert+\kappa(x)\vert s_0\vert,$$
where $\kappa$ is defined in {\bf A1}.
We next discuss {\bf A3}.
\begin{enumerate}
\item
For the probit model, we have $F(\lambda)=\Phi(\lambda)=\int_{-\infty}^{\lambda}\frac{1}{\sqrt{2\pi}} e^{-u^2/2}\,du$.
We have that
\begin{eqnarray*}
d_{TV}(p(\cdot \mid \lambda), p(\cdot \mid \lambda'))
& = & \frac{1}{2} \; \left\{ |(1-\Phi(\lambda))-(1-\Phi(\lambda'))| \,+\, |\Phi(\lambda)-\Phi(\lambda')| \right\} \\
& = & |\Phi(\lambda)-\Phi(\lambda')|.
\end{eqnarray*}
In this case, we have for any $x\in \R$ and $h\geq 0$, $\left\vert \Phi(x+h)-\Phi(x)\right\vert\leq 2\Phi(h/2)-1=1-2(1-\Phi(h/2))$. 
Around $h=0$, we have $2(1-\Phi(h/2))\geq \exp(-d_1 h)$, provided that $d_1>1/\sqrt{2\pi}$.
Fix such $d_1$ and suppose that such the latter inequality is valid for $0\leq h\leq \epsilon=\epsilon(d_1)$.
Next using an integration by parts, we have for $h\geq 2$,
$$2\left(1-\Phi(h/2)\right)\geq \frac{2\exp(-h^2/8)}{\sqrt{2\pi}h}.$$
If $d_2>1/8$ is fixed, one can choose $M>\epsilon$ such that for $h\geq M$,  
$$2\left(1-\Phi(h/2)\right)\geq \exp(-d_2 h^2).$$
For $\epsilon\leq h\leq M$, we have for a sufficiently large real number $d_3$,
$$2\left(1-\Phi(h/2)\right)\geq 2\left(1-\Phi(M/2)\right)\geq \exp(-d_3\epsilon^2)\geq \exp(-d_3 h^2).$$  
Setting $d=\max(d_1,d_2,d_3)$, one can choose $\phi(h)=d(h+h^2)$.
\item
For the logistic model, we have $F(\lambda)=\exp(\lambda)/\left(1+\exp(\lambda)\right)$. 
In this case, one can use $\phi(s)=\vert s\vert$ in {\bf A3}. A proof will be given below directly for the multinomial case.$\square$
\end{enumerate}

Now we extend our results to categorical time series.
Set $E=\{0,1,\ldots,N-1\}$ and $F=\R^{N-1}$. For $i=1,\ldots, N-1$, we set $p(i\vert s)=\frac{\exp(s_i)}{S(s)}$ with $S(s)=1+\sum_{j=1}^{N-1}\exp(s_j)$.
Then $p(0\vert s)=S(s)^{-1}$. The corresponding model, called multinomial autoregressive model, is introduced for instance in \citep{russell}.

\begin{prop}\label{categorical}
Suppose that Assumption {\bf A1} holds true and there exists $s_0\in F$ such that $x\mapsto f(s_0,y,x)$ is in $\mathcal{M}_{\log}$ for any $y\in E$. Then the conclusions of Theorem \ref{main} are valid 
for multinomial autoregressive model.
\end{prop}

\paragraph{Proof of Proposition \ref{categorical}}
Checking {\bf A2} is exactly as in the proof of Proposition \ref{binary}. We then check {\bf A3}. Observe that
$$d_{TV}\left(p(\cdot\vert s),p(\cdot\vert s')\right)=\frac{1}{2}\sum_{i=0}^{N-1}\left\vert p(i\vert s)-p(i\vert s')\right\vert.$$
w.l.o.g. we assume that $S(s)\geq S(s')$ and we set
$$I_+=\left\{1\leq i\leq N-1: p(i\vert s)>p(i\vert s')\right\}$$
and $I_{-}=\{1,\ldots, N-1\}\setminus I_+$.
We have 
\begin{eqnarray*}
d_{TV}\left(p(\cdot\vert s),p(\cdot\vert s')\right)&=&\frac{1}{2}\sum_{i\in I_+}\left[p(i\vert s)-p(i\vert s')\right]+\frac{1}{2}\sum_{i\in I_-}\left[p(i\vert s')-p(i\vert s)\right]\\
&+& S(s')^{-1}-S(s)^{-1}\\
&=&\sum_{i\in I_+}\left[p(i\vert s)-p(i\vert s')\right]\\
&=& \sum_{i\in I_+}p(i\vert s)\left[1-\frac{S(s)}{S(s')}e^{s'_i-s_i}\right]\\
&\leq& \sum_{i\in I_+}p(i\vert s)\left[1-e^{s'_i-s_i}\right]\\
&\leq& \max_{i\in I_+}\left[1-e^{s'_i-s_i}\right].
\end{eqnarray*}
When $s\neq s'$, $I_+$ is not empty and if $i\in I_+$, then $e^{s_i-s'_i}>\frac{S(s)}{S(s')}\geq 1$ and then $s_i>s'_i$. 
We conclude that 
$$\max_{i\in I_+}\left[1-e^{s'_i-s_i}\right]\leq 1-e^{-\Vert s-s'\Vert_{\infty}}.$$
This proves {\bf A3} with $\phi(h)=h$ and the infinite norm on $L$ (and then any norm on $L$ by equivalence).$\square$

\subsection{Count time series}

We first consider the Poisson conditional distribution, $p(\cdot\vert s)=\mbox{Pois}(s)$. Here we set $E=\N$ and $F=\R_+$.
We then have
$$\P\left(Y_t\in A\vert X,Y_{t-1},\lambda_t,Y_{t-2},\lambda_{t-1}\ldots\right)=\mbox{Pois}\left(\lambda_t\right),\quad \lambda_t=f\left(\lambda_{t-1},Y_{t-1},X_{t-1}\right).$$

Note that here, $\int yp(dy\vert s)=s$ and one can check the assumptions of Proposition \ref{pratique1} or Proposition \ref{pratique2} with $V(s)=1+s$.

\begin{prop}\label{Poisson}
Suppose that Assumptions {\bf A1} and {\bf A2(1)} hold true. Then the conclusions of Theorem \ref{main} are valid.
\end{prop}

\paragraph{Note.} 
When, the function $f$ is linear, i.e. $f(s,y,x)=\kappa(x)s+\widetilde{\kappa}(x)y+\gamma(x)$, with nonnegative functions $\kappa,\widetilde{\kappa},\gamma$, we get an extension of the classical INGARCH model considered in \citep{FTR} by allowing exogenous covariates in the dynamic. But using Proposition \ref{pratique2}, one can also deal with Poisson threshold autoregressive processes.
This kind of model has been considered without exogenous covariate for instance in \citep{DDM} or \citep{Wang}. In this case, 
our assumptions are similar for getting existence of an ergodic solution.  We then also get a non-trivial extension by allowing exogenous covariates in the random intensity $\lambda_t$ of this model.

\paragraph{Proof of Proposition \ref{Poisson}}
From Proposition \ref{pratique1}, Assumption {\bf A2} is satisfied. We only need to check {\bf A3}. For $s<s'$ and independent $U\sim \mbox{Pois}(s)$,
$V\sim \mbox{Pois}(s'-s)$, we have that $U+V\sim \mbox{Pois}(s')$. Therefore,
\begin{eqnarray*}
d_{TV}(p(\cdot \vert s), p(\cdot \vert s'))
& \leq & P( V\neq 0) \\
& = & 1 \,-\, e^{-|s-s'|}.
\end{eqnarray*}
Hence {\bf A3} is satisfied with $\phi(h)=h$.$\square$
\bigskip

Next, we study another count autoregressive model which is quite popular because it replaces the Poisson distribution by a distribution that takes into account the over-dispersion of count data. This model, called negative binomial, is studied for instance in \citep{Davis}. We remind that the negative binomial distribution NB$(r,q)$ with parameter $r\in \N^{*}$ and $q\in (0,1)$ 
can be defined as a mixture of Poisson distribution, for instance it equals the probability distribution of the random variable $N_{\frac{\varepsilon}{r}s}$ where $\varepsilon$ follows a gamma distribution with parameters $(r,1)$ and is independent from a Poisson process $N$ with intensity $1$ and $s=\frac{qr}{1-q}$, which equals to the mean of this distribution.
Here for a given positive integer $r$, we assume that
$$p\left(\cdot\vert s\right)=\mbox{NB}\left(r,\frac{s}{s+r}\right).$$

\begin{prop}\label{nebin}
Let Assumptions {\bf A1-A2(1)} hold true. Then the conclusions of Theorem \ref{main} are valid for the negative binomial autoregressive model.
\end{prop}

\paragraph{Proof of Proposition \ref{nebin}}
Since $\int y p(dy\vert s)=s$, Proposition \ref{pratique1} ensures the validity of {\bf A2}.
We then check {\bf A3}. Denoting by $f_{\varepsilon}$ the probability density of $\varepsilon$, we have
\begin{eqnarray*}
d_{TV}\left(NB(r,q),NB(r,q')\right)&\leq & \P\left(N_{\frac{\varepsilon}{r}s}\neq N_{\frac{\varepsilon}{r}s'}\right)\\
&=& 1-\P\left(N_{\frac{\varepsilon}{r}\vert s-s'\vert}=0\right)\\
&=& 1-\int \exp\left(-\frac{u}{r}\vert s-s'\vert\right)f_{\varepsilon}(u)du\\
&=& 1- \left(\frac{1}{1+\frac{1}{r}\vert s-s'\vert}\right)^r\\
&\leq & 1-\exp\left(-\vert s-s'\vert\right),
\end{eqnarray*}
using the expression of the Laplace transform of the gamma distribution. We then get {\bf A3} with $\phi(h)=h$ and the proof of the proposition is then complete.$\square$

\subsection{GARCH type processes}
GARCH processes are defined by the recursions
$$Y_t=\varepsilon_t\sigma_t,\quad \lambda_t:=\sigma^2_t=f\left(\sigma^2_{t-1},Y_{t-1}X_{t-1}\right),$$
where $\left(\varepsilon_t\right)_{t\in\Z}$ a sequence of i.i.d. random variables such that $\E\left(\varepsilon_0\right)=0$, $\E\left(\varepsilon_0^2\right)=1$ and the sequences $\left(\varepsilon_t\right)_{t\in \Z}$ and $\left(X_t\right)_{t\in \Z}$ are independent.

\begin{prop}\label{GARCH}
Suppose that Assumptions {\bf A1-A2(2)} hold true with $f$ lower bounded by a positive constant $c_{-}$. Assume furthermore that the noise $\varepsilon_0$ has a probability density $f_{\varepsilon}$ non-decreasing
on $(-\infty,0]$ and non-increasing on $[0,\infty)$. Then the conclusions of Theorem \ref{main} hold true.
\end{prop}

\paragraph{Notes}
\begin{enumerate}
\item
The main restriction with our approach is the additional constraint on the variation of the density $f_{\varepsilon}$. Without this restriction, we did not find an argument for checking {\bf A3}. 
This condition on the density, which is probably not optimal, is not so restrictive.  One can always consider standard symmetric densities such as that of Gaussian, Laplace or Student distributions. But non-symmetric densities are also possible.  
\item
Despite our restriction on the noise density, our approach can be used to define models with a complex structure for the conditional variance $\sigma^2_t$, in particular threshold models. There exist several versions of ARCH or GARCH threshold models in the literature. See for instance \citep{Gourieroux}, \citep{Cline} or \citep{Zak} for ARCH versions and \citep{Zak2} for a GARCH version.
Our version allows exogenous covariates and the threshold is not necessarily $0$ as in \citep{Zak2}. Moreover, it is not difficult to generalize Proposition \ref{pratique2} to allow multiple threshold, in the spirit of \citep{Gourieroux} for the ARCH. In this case, our result is also interesting even without exogenous covariates.
\end{enumerate}

\paragraph{Proof of Proposition \ref{GARCH}}
From Proposition \ref{pratique1}, Assumption {\bf A2} holds true. 
It is only necessary to check {\bf A3}. Here $F=[c_{-},\infty)$ and we have $p(dy\vert s)=\frac{1}{\sqrt{s}}f_{\varepsilon}\left(y/\sqrt{s}\right)dy$.
Since $f_{\varepsilon}$ is non-decreasing on $(-\infty,0)$ and non-increasing on $(0,\infty)$, we have 
for $c_{-}^2\leq s'\leq s$ and $w=\sqrt{s}$, $w'=\sqrt{s'}$, 
\begin{eqnarray*}
\frac{1}{2}\int \left\vert\frac{1}{w}f_{\varepsilon}\left(\frac{u}{w}\right)-\frac{1}{w'}f_{\varepsilon}\left(\frac{u}{w'}\right)\right\vert du
&=& 1-\int \left(\frac{1}{x}f_{\varepsilon}\left(\frac{u}{w}\right)\right)\wedge \left(\frac{1}{w'}f_{\varepsilon}\left(\frac{u}{w'}\right)\right)du\\
&\leq &1-\frac{1}{w}\int f_{\varepsilon}\left(\frac{u}{w}\right)\wedge f_{\varepsilon}\left(\frac{u}{w'}\right)du\\  
&=&1-\frac{1}{w}\int f_{\varepsilon}\left(\frac{u}{w'}\right)du\\
&=& 1-\frac{w'}{w}\\
&\leq& 1-\exp\left(-\frac{w-w'}{c_{-}}\right)\\
&\leq& 1-\exp\left(-\frac{s-s'}{2c_{-}^{3/2}}\right).
\end{eqnarray*}
This completes the proof.$\square$

\subsection{Conditionally homoscedastic autoregressive processes}
Finally we consider the transition kernel  
$$p\left(dy\vert s\right)=f_{\varepsilon}(y-s)dy,$$
where $f_{\varepsilon}$ denotes the probability density of a random variable $\varepsilon_0$.
This case covers the model
\begin{equation}\label{end}
Y_t=\lambda_t+\varepsilon_t,\quad \lambda_t=f\left(\lambda_{t-1},Y_{t-1},X_{t-1}\right)
\end{equation}
when the two processes $\left(X_t\right)_{t\in\Z}$ and $\left(\varepsilon_t\right)_{t\in\Z}$ are independent and $\left(\varepsilon_t\right)_{t\in\Z}$ is a sequence of i.i.d. random variables with probability density $f_{\varepsilon}$.
Note that when $f(s,y,x)=g(y,x)$, we obtain an autoregressive model that includes for instance the well-known threshold autoregressive model, see for instance \citep{Tong} and \citep{Tsay}, with lag $1$ but with exogenous covariates.
When $f(s,y,x)=a(x)s+g(y,x)-a(x)y$, we obtain
$$Y_t=g\left(Y_{t-1},X_{t-1}\right)+\varepsilon_t-a\left(X_{t-1}\right)\varepsilon_{t-1},$$
which includes ARMA(1,1) models with varying-coefficients.

\begin{prop}\label{end}
Suppose that Assumptions {\bf A1-A2(1)} hold true, $\E\vert \varepsilon_0\vert<\infty$, the density $f_{\varepsilon}$ is symmetric around $0$, continuous at point $0$, non-increasing on $(0,\infty)$ and non-decreasing on $(-\infty ,0]$ and such that for some $h_1,C,C'>0$ and a positive integer $C''>2$ such that
$$f(y)\geq C \exp\left(-C'y^{C''}\right),\quad y\geq h_1>0.$$
Then the conclusions of Theorem \ref{main} are valid.
\end{prop}

\paragraph{Proof of Proposition \ref{end}}
Using Proposition \ref{pratique1}, one can check {\bf A2}.  
To check Assumption {\bf A3}, we use the equalities
$$d_{TV}\left(p(\cdot\vert s),p(\cdot\vert s')\right)=1-\int f_{\varepsilon}(y-s)\wedge f_{\varepsilon}(y-s')dy=1-\int f_{\varepsilon}(y)\wedge f_{\varepsilon}(y-(s'-s))dy.$$
We then need to derive a lower bound for
$$I(h)=\int f_{\varepsilon}(y)\wedge f_{\varepsilon}(y-h)dy,\quad h\geq 0.$$
Using the symmetry of the density, it is easily seen that 
$$I(h)=\int_{h/2}^{\infty}f_{\varepsilon}(y)dy+\int_{-\infty}^{h/2}f_{\varepsilon}(y-h)dy=2-2F(h/2),$$
where $F(h)=\int_{-\infty}^h f_{\varepsilon}(y)dy$.
Note next that for $h\geq h_1$, we have 
$$1-F(h)\geq J(h)=\int_{h}^{\infty}C\exp\left(-C'y^{C''}\right)dy=\int_{h}^{\infty}Cy^{{C''}-1}\exp\left(-C'y^{C''}\right)y^{1-C''}dy$$
and using an integration by part, we have
$$J(h)\geq \frac{C\exp\left(-C'h^{C''}\right)}{C'{C''}h^{{C''}-1}}- \frac{J(h)}{C'{C''}\left({C''}-2\right)h_1^{{C''}-2}}.$$
Then there exists a positive constant $D'>0$ large enough such that
$$J(h)\geq \frac{1}{2}\exp\left(-C'h^{D'}\right).$$
Moreover, at point $0$, we have 
$1-F(h)=1/2-f_{\varepsilon}(0)h+o(h)$.
If $C_2>2f_{\varepsilon}(0)$, there then exists $h_0\in (0,h_1)$  
$$1-F(h)\geq \frac{1}{2}\exp\left(-C_2 h\right),\quad 0\leq h\leq h_0.$$
Finally, if $h_0\leq h\leq h_1$, we have 
$$1-F(h)\geq 1-F(h_1)\geq \frac{1}{2}\exp\left(-C' h_1^{D'}\right)\geq \frac{1}{2}\exp\left(-C'\frac{h_1^{D'}}{h_0^{D'}}h^{D'}\right).$$
Setting $D=\max\left(C_2/2,C' h_1^{D'}/(2h_0)^{D'}\right)$, we get 
$$I(h)\geq \exp\left(-D\left(h+h^{D'}\right)\right),\quad h\geq 0,$$
which gives {\bf A3} with $\phi(h)=D\left(h+h^{D'}\right)$. $\square$

\bibliographystyle{plainnat}
\bibliography{biblogistic}

\end{document}